\documentclass[11pt,a4paper]{article}
\usepackage{a4wide,amsfonts,amsmath,latexsym,amssymb,euscript,eufrak,graphicx,units,mathrsfs}

\usepackage{float}
\newfloat{figure}{H}{lof}
\floatname{figure}{\figurename}

\DeclareMathAlphabet{\eufrak}{U}{}{}{}
\SetMathAlphabet\eufrak{normal}{U}{euf}{m}{n}
\SetMathAlphabet\eufrak{bold}{U}{euf}{b}{n}

\usepackage{amsmath, amsthm, amsfonts, amssymb}

\numberwithin{equation}{section}

\def\real{{\mathord{{\rm I\kern-2.8pt R}}}}
\def\inte{{\mathord{{\rm I\kern-2.8pt N}}}}
\def\PP{{\mathord{{\rm I\kern-2.8pt P}}}}

\def\real{{\mathord{\mathbb R}}}

\def\inte{{\mathord{\mathbb N}}}

\def\Dom{{\mathrm{{\rm Dom}}}}

\def\HH{\EuFrak H}

\newcommand{\e}{\varepsilon}

\def\E{\mathop{\hbox{\rm I\kern-0.20em E}}\nolimits}

\def\real{\mathbb{R}}

\newtheorem{prop}{Proposition}[section]

\newtheorem{lemma}[prop]{Lemma}
\newtheorem{definition}[prop]{Definition}
\newtheorem{corollary}[prop]{Corollary}
\newtheorem{theorem}[prop]{Theorem}
\newtheorem{thm}[prop]{Theorem}
\newtheorem{remark}[prop]{Remark}

\allowdisplaybreaks

\begin{document}

\begin{center}
{\large \textbf{Multivariate normal approximation using Stein's method and Malliavin calculus}}\\[0pt]
~\\[0pt]
Ivan Nourdin\footnote{LPMA,
Universit\'e Paris VI, Bo\^ite courrier
188, 4 place Jussieu, 75252
Paris Cedex 05, France. Email: \texttt{ivan.nourdin@upmc.fr}},
Giovanni Peccati\footnote{
\'Equipe Modal'X, Universit\'e Paris Ouest - Nanterre la D\'efense, 200 avenue de la R\'epublique, 92000 Nanterre
and LSTA, Universit\'e Paris VI, France. Email: \texttt{giovanni.peccati@gmail.com}}
and Anthony R\'eveillac\footnote{
Institut f\"ur Mathematik, Humboldt-Universit\"at zu Berlin,
Unter den Linden 6, 10099 Berlin, Germany. Email:
\texttt{areveill@mathematik.hu-berlin.de}}\\[0pt]
{\it Universit\'e Paris VI, Universit\'e Paris Ouest and Humboldt-Universit\"at zu Berlin}\\~\\
~\\[0pt]
\end{center}

{\small \noindent \textbf{Abstract:} We combine Stein's method
with Malliavin calculus in order to obtain explicit bounds in the
multidimensional normal approximation (in the Wasserstein
distance) of functionals of Gaussian fields. Our results
generalize and refine the main findings by Peccati and Tudor
(2005), Nualart and Ortiz-Latorre (2007), Peccati (2007) and
Nourdin and Peccati (2007b, 2008); in particular, they apply to
approximations by means of Gaussian vectors with an arbitrary,
positive definite covariance matrix. Among several examples, we
provide an application to a functional version of the Breuer-Major
CLT for fields subordinated to a fractional Brownian motion.\\

\normalsize
}

{\small \noindent \textit{Key words:} Breuer-Major CLT, fractional
Brownian motion, Gaussian processes, Malliavin calculus, Normal
approximation, Stein's method, Wasserstein distance.~\newline
\normalsize }

{\small \noindent \textbf{R\'esum\'e:} Nous expliquons comment combiner la m\'ethode
de Stein avec les outils du calcul de Malliavin
pour majorer, de mani\`ere explicite, la distance de Wasserstein
entre une fonctionnelle
d'un champs gaussien donn\'ee et son approximation normale multidimensionnelle.
Notre travail g\'en\'eralise et affine des r\'esultats ant\'erieurs prouv\'es par
Peccati et Tudor
(2005), Nualart et Ortiz-Latorre (2007), Peccati (2007) et
Nourdin et Peccati (2007b, 2008).
Entre autres exemples, nous associons des bornes \`a la version
fonctionnelle du th\'eor\`eme de la limite centrale de Breuer-Major dans le cas du
mouvement brownien fractionnaire.\\
\normalsize
}

{\small \noindent \textit{Mots cl\'es:} th\'eor\`eme de la limite centrale de Breuer et Major, 
mouvement brownien fractionnaire, processus gaussiens, calcul de Malliavin,
approximation normale, m\'ethode de Stein, distance de Wasserstein.~\newline
\normalsize }

{\small \noindent \textit{This revised version}: November 2008}
\normalsize

\section{Introduction}
Let $Z\sim\mathscr{N}(0,1)$ be a standard Gaussian random variable
on some probability space $(\Omega,\mathcal{F},P)$, and let $F$ be
a real-valued functional of an infinite-dimensional Gaussian
field. In the papers \cite{NourdinPeccati1, NourdinPeccati2} it is
shown that one can combine Stein's method (see e.g.
\cite{Chen_Shao_sur}, \cite{Stein1} or \cite{Stein2}) with
Malliavin calculus (see e.g. \cite{Nualart}), in order to deduce
explicit (and, sometimes, optimal) bounds for quantities of the
type $d(F,Z)$, where $d$ stands for some distance between the law
of $F$ and the law of $Z$ (e.g., $d$ can be the Kolmogorov or the
Wasserstein distance). The aim of this paper is to extend the
results of \cite{NourdinPeccati1, NourdinPeccati2} to the
framework of the {\sl multidimensional} Gaussian approximation in
the Wasserstein distance. Once again, our techniques hinge upon
the use of infinite-dimensional operators on Gaussian spaces (like
the \textsl{divergence operator} or the \textsl{Ornstein-Uhlenbeck
generator}) and upon an appropriate multidimensional version of
Stein's method (in a form close to Chatterjee and Meckes
\cite{ChatterjeeMeckes}, but see also Reinert and R\"{o}llin
\cite{ReinertRollin}). As a result, we will obtain explicit
bounds, both in terms of Malliavin derivatives and contraction
operators, thus providing a substantial refinement of the main
findings by Nualart and Ortiz-Latorre \cite{NualartOrtizLatorre}
and Peccati and Tudor \cite{PeccatiTudor}. Note that an important
part of our computations (see e.g. Lemma \ref{lm-control}) are
directly inspired by those contained in
\cite{NualartOrtizLatorre}: we shall indeed stress that this last
reference contains a fundamental methodological breakthrough,
showing that one can deal with (possibly multidimensional) weak
convergence on a Gaussian space, by means of Malliavin-type
operators and ``characterizing'' differential equations. See
\cite{NouPec07} for an application of these techniques to
non-central limit theorems. Incidentally, observe that the paper
\cite{Pecc}, which is mainly based on martingale-type techniques,
also uses distances between probability measures (such as the
Prokhorov distance) to deal with multidimensional Gaussian
approximations on Wiener space, but without giving explicit
bounds.

The rationale behind Stein's method is better understood in
dimension one. In this framework, the starting point is the
following crucial result, proved e.g. in \cite{Stein1}.
\begin{lemma}[Stein's Lemma]
\label{lemma:SteinIBP1}
A random variable $Y$ is such that $Y
\overset{\rm Law}{=}Z\sim \mathscr{N}(0,1)$ if and only if, for
every continuous and piecewise continuously differentiable
function $f:\mathbb{R}\to\mathbb{R}$
such that $E\big|f'(Z)\big|<\infty$, one has
\begin{equation}\label{urStein}
 E[f'(Y)-Yf(Y)]=0.
\end{equation}
\end{lemma}
The fact that a random variable $Y$ satisfying (\ref{urStein}) is
necessarily Gaussian can be proved by several routes: for
instance, by taking $f$ to be a complex exponential, one can show
that the characteristic function of $Y$, say $\psi(t)$, is
necessarily a solution to the differential equation
$\psi'(t)+t\psi(t)=0$, and therefore $\psi(t)=\exp(-t^2/2)$;
alternatively, one can set $f(x)=x^n$, $n=1,2,...$, and observe
that (\ref{urStein}) implies that, for every $n$, one must have
$E(Y^n)=E(Z^n)$, where $Z\sim \mathscr{N}(0,1)$ (note that the law
of $Z$ is determined by its moments).

Heuristically, Lemma \ref{lemma:SteinIBP1} suggests that the
distance $d(Y,Z)$, between the law of a random variable $Y$ and
that of $Z\sim \mathscr{N}(0,1)$, must be ``small'' whenever
$ E[f'(Y)-Yf(Y)] \simeq 0, $
for a sufficiently large class of functions $f$. In the seminal
works \cite{Stein1, Stein2}, Stein proved that this somewhat
imprecise argument can be made rigorous by means of the use of
differential equations. To see this, for a given function
$g:\real\to\real$, define the \textsl{Stein equation} associated
with $g$ as
\begin{equation}
\label{eq:SteinEquation1}
g(x)-E[g(Z)] = h'(x)-x h(x), \quad \forall x\in \real,
\end{equation}
(we recall that $Z\sim\mathscr{N}(0,1)$). A solution to
(\ref{eq:SteinEquation1}) is a function $h$ which is
Lebesgue-almost everywhere differentiable, and such that there
exists a version of $h'$ satisfying (\ref{eq:SteinEquation1}) for
every $x\in\real$.  If one assumes that $g\in {\rm Lip}(1)$ (that
is, if $\|g\|_{Lip}\leq 1$, where $\|\cdot\|_{Lip}$ stands for the
usual Lipschitz seminorm), then a standard result (see e.g.
\cite{Stein2}) yields that (\ref{eq:SteinEquation1}) admits a
solution $h$ such that $\|h'\|_\infty \leq 1$ and $\|h''\|_\infty
\leq 2$. Now recall that the \textsl{Wasserstein distance} between
the laws of two real-valued random variables $Y$ and $X$ is
defined as
\begin{equation*}
d_{\rm W}(Y,X) = \sup_{g\in {\rm Lip}(1)}
\left\vert E[g(Y)]-E[g(X)]\right\vert,
\end{equation*}
and introduce the notation $\mathscr{F}_{\rm W} = \{ f:
\|f'\|_\infty \leq 1, \, \|f''\|_\infty \leq 2\}$. By taking
expectations on the two sides of (\ref{eq:SteinEquation1}), one
obtains finally that, for $Z\sim\mathscr{N}(0,1)$ and for a
generic random variable $Y$,
\begin{equation}
\label{Stbound}
d_{\rm W}(Y,Z) \leq \sup_{f\in \mathscr{F}_{\rm W}} \left\vert
E[f'(Y)-Yf(Y)]\right\vert,
\end{equation}
thus giving a precise meaning to the heuristic argument sketched
above (note that an analogous conclusion can be obtained for other
distances, such as the total variation distance or the Kolmogorov
distance -- see e.g. \cite{Chen_Shao_sur} for a discussion of this
point). We stress that the topology induced by $d_{\rm W}$, on
probability measures on $\real$, is stronger than the topology
induced by weak convergence.

The starting point of \cite{NourdinPeccati1,NourdinPeccati2} is
that a relation such as (\ref{Stbound}) can be very effectively
combined with Malliavin calculus, whenever $Y$ is a centered
regular functional of some infinite dimensional Gaussian field.
To see this, denote by $DY$ the Malliavin derivative of $Y$
(observe that $DY$ is a random element with values in some
adequate Hilbert space $\EuFrak{H}$), and write $L$ to indicate
the (infinite-dimensional) Ornstein-Uhlenbeck generator (see
Section \ref{section:Malliavin} below for precise definitions).
One crucial relation proved in \cite{NourdinPeccati1}, and then
further exploited in \cite{NourdinPeccati2}, is the upper bound
\begin{equation}
\label{GioIvan}
d_{\rm W}(Y,Z) \leq E|1-\langle DY,
-DL^{-1}Y\rangle_\EuFrak{H}|.
\end{equation}

As shown in \cite{NourdinPeccati1}, when specialized to the case
of $Y$ being equal to a multiple Wiener-It\^{o} integral,
relation (\ref{GioIvan}) yields bounds that are intimately related
with the CLTs proved in \cite{NualartOrtizLatorre} and \cite{NP}.
See \cite{NourdinPeccati2} for a characterization of the
optimality of these bounds; see again \cite{NourdinPeccati1} for
extensions to non-Gaussian approximations and for applications to
the Breuer-Major CLT (stated and proved in \cite{BM}) for
functionals of a fractional Brownian motion.

The principal contribution of the present paper (see e.g. the
statement of Theorem \ref{theo:majDist} below) consists in
showing that a relation similar to (\ref{GioIvan}) continues to
hold when $Z$ is replaced by a
$d$-dimensional ($d\geq2$) Gaussian vector
$F=(F_1,...,F_d)$ of smooth functionals of
a Gaussian field, and $d_{\rm W}$ is the Wasserstein distance
between probability laws on $\real^d$ (see Definition
\ref{DefWass} below). Our results apply to Gaussian approximations
by means of Gaussian vectors with arbitrary positive definite
covariance matrices. The proofs rely on a multidimensional version
of the Stein equation (\ref{eq:SteinEquation1}), that we combine
with standard integration by parts formulae on an
infinite-dimensional Gaussian space. Our approach bears some
connections with the paper by Hsu \cite{Hsu}, where the author
proves an hybrid Stein/semimartingale characterization of Brownian
motions on manifolds, via Malliavin-type operators.

The paper is organized as follows. In Section
\ref{section:Malliavin} we provide some preliminaries on Malliavin
calculus. Section \ref{S : Main} contains our main results,
concerning Gaussian approximations by means of vectors of Gaussian
random variables with positive definite covariance matrices.
Finally, Section \ref{section:applications} deals with two
applications: (i) to a functional version of the Breuer-Major CLT
(see \cite{BM}), and (ii) to Gaussian approximations of
functionals of finite normal vectors, providing a generalization
of a technical result proved by Chatterjee in
\cite{Chatterjee_ptrf}.

\section{Preliminaries and notation}
\label{section:Malliavin}
In this section, we recall some basic
elements of Malliavin calculus for Gaussian processes. The reader
is referred to \cite{Nualart} for a complete discussion of this
subject. Let $X=\{X(h),\; h\in \EuFrak{H}\}$ be an
\textsl{isonormal Gaussian process} on a probability space
$(\Omega,\mathcal{F},P)$. This means that $X$ is a centered
Gaussian family indexed by the elements of an Hilbert space
$\EuFrak{H}$, such that, for every pair $h,g\in\EuFrak{H}$ one has
that
$ E[X(h) X(g)]=\langle h, g\rangle_{\EuFrak{H}}.$

We let $L^2(X)$ be shorthand for the space
$L^2(\Omega,\sigma(X),P)$. It is well known that every random
variable $F \in L^2(X)$ admits the chaotic expansion
$ F=E(F)+\sum_{n=1}^\infty I_n(f_n) $
where the deterministic kernels $f_n$, $n\geq 1$, belong to
$\EuFrak{H}^{\odot n}$ and the convergence of the series holds in
$L^2(X)$. One sometimes uses the notation $I_0(f_0)=E[F]$. In the
particular case where $\EuFrak{H}:=L^2(T,\mathcal{A},\mu)$, with
$(T,\mathcal{A})$ a measurable space
and $\mu$ is a $\sigma$-finite measure without atoms,
the random variable
$I_n(f_n)$ coincides with the \textsl{multiple Wiener-It\^o
integral} (of order $n$) of $f_n$ with respect to $X$ (see \cite[Section 1.1.2.]{Nualart}).

Let $f\in\EuFrak{H}^{\odot p}$, $g\in\EuFrak{H}^{\odot q}$ and
$0 \leq r \leq p\wedge q$. We define the $r$th \textsl{contraction}
$f\otimes_r g$
of $f$ and $g$ as the element
of $\EuFrak{H}^{\otimes (p+q-2r)}$ given by
$$ f\otimes_r g :=\sum_{i_1,\ldots,i_r=1}^\infty \langle f, e_{i_1}\otimes
\ldots\otimes e_{i_r} \rangle_{\EuFrak{H}^{\otimes r}} \otimes
\langle g, e_{i_1} \otimes \ldots\otimes e_{i_r}
\rangle_{\EuFrak{H}^{\otimes r}},$$ where $\{e_k,\; k\geq 1\}$ is
a complete orthonormal system in $\EuFrak{H}$. Note that
$f\otimes_0 g=f\otimes g$; also, if $p=q$, then
$f\otimes_p g=\langle f, g \rangle_{\EuFrak{H}^{\otimes
p}}$. Note that, in general, $f\otimes_r g$ is not a symmetric
element of $\EuFrak{H}^{\otimes (p+q-2r)}$; the canonical
symmetrization of $f\otimes_r g$ is denoted by $f
\widetilde{\otimes}_r g$. We recall the
 product formula for multiple stochastic integrals:

\begin{equation*}
I_p(f)I_q(g)=\sum_{r=0}^{p\wedge q}
r! \binom{p}{r}\binom{p}{q}
I_{p+q-2r}(f\widetilde{\otimes}_r g).
\end{equation*}

Now, let $\mathscr{S}$ be the set of cylindrical functionals $F$ of
the form
\begin{equation}
\label{eq:cylindrical}
F=\varphi(X(h_1),\ldots,X(h_n)),
\end{equation}
where $n\geq 1$, $h_i \in \EuFrak{H}$ and the function $\varphi\in
\mathscr{C}^{\infty}(\real^n)$ is such that its partial
derivatives have polynomial growth. The \textsl{Malliavin
derivative} $DF$ of a functional $F$ of the form
(\ref{eq:cylindrical}) is the square integrable
$\EuFrak{H}$-valued random variable defined as
\begin{equation}
\label{baldassarre}
DF=\sum_{i=1}^n \partial_i \varphi(X(h_1),\ldots,X(h_n)) h_i,
\end{equation}
where $\partial_i \varphi$ denotes the $i$th partial derivative of $\varphi$.
In particular, one has that $DX(h)=h$ for every $h$ in
$\EuFrak{H}$. By iteration, one can define the $m$th derivative
$D^m F$ of $F\in\mathscr{S}$, which is an element of
$L^2(\Omega;\EuFrak{H}^{\odot m})$,
for $m\geq 2$. As usual
$\mathbb{D}^{m,2}$ denotes the closure of $\mathscr{S}$ with
respect to the norm $\|\cdot\|_{m,2}$ defined by the relation
$ \|F\|_{m,2}^2=E[F^2] +\sum_{i=1}^m E[\|D^iF\|^2_{\EuFrak{H}^{\otimes i}}].$

Note that every finite sum of Wiener-It\^o integrals
always belongs to $\mathbb{D}^{m,2}$ ($\forall m\geq 1$).
The Malliavin derivative $D$ satisfies the following
\textsl{chain rule formula}: if $\varphi:\real^n\to \real$ is in
$\mathscr{C}_b^1$ (defined as the set of continuously
differentiable functions with bounded partial derivatives) and if
$(F_1,\ldots,F_n)$ is a random vector such that each component
belongs to $\mathbb{D}^{1,2}$, then $\varphi(F_1,\ldots,F_n)$ is
itself an element of $\mathbb{D}^{1,2}$, and moreover
\begin{equation}
\label{eq:ChainRule}
D\varphi(F_1,\ldots,F_n)=\sum_{i=1}^n \partial_i\varphi(F_1,\ldots,F_n) DF_i.
\end{equation}

The \textsl{divergence operator} $\delta$ is defined as the dual
operator of $D$. Precisely, a random element $u$ of $L^2(\Omega;\EuFrak{H})$
belongs to the domain of $\delta$ (denoted by ${\rm Dom} \delta$)
if there exists a constant $c_u$ satisfying
$\left\vert E[\langle DF, u \rangle_{\EuFrak{H}}]\right\vert \leq c_u \|F\|_{L^2(\Omega)}$
for every  $F \in \mathscr{S};$
in this case, the divergence of $u$, written $\delta(u)$, is
defined by the following duality property:
\begin{equation}
\label{eq:MalliavinIBP}
E[F \delta(u)]=E[\langle DF, u
\rangle_{\EuFrak{H}}], \quad \forall F \in \mathbb{D}^{1,2}.
\end{equation}
The crucial relation (\ref{eq:MalliavinIBP}) is customarily called
the (Malliavin) \textsl{integration by parts formula}.

In what follows, we shall denote by $T = \{T_t: t\geq 0\}$ the
\textsl{Ornstein-Uhlenbeck semigroup}. We recall that, for every
$t\geq 0$ and every $F\in L^2(X)$,
\begin{equation}
\label{OrnstU}
T_t(F) = \sum_{n=0}^\infty e^{-nt}J_n(F),
\end{equation}
where, for every $n\geq0$ and for the rest of the paper, the
symbol $J_n$ denotes the projection operator onto the $n$th Wiener
chaos,
that is onto the closed linear subspace of $L^2(X)$ generated by the random variables
of the form $H_n(X(h))$ with $h\in\EuFrak{H}$ such that $\|h\|_\EuFrak{H}=1$,
and $H_n$ the $n$th Hermite polynomial defined by (\ref{her-pol}).
Note that $T$ is indeed the semigroup associated with an
infinite-dimensional stationary Gaussian process with values in
$\real^{\HH}$, having the law of $X$ as an invariant distribution
(see e.g. \cite[Section 1.4]{Nualart} for a more detailed
discussion of the Ornstein-Uhlenbeck semigroup in the context of
Malliavin calculus; see Barbour \cite{Barbour} for a version of
Stein's method involving Ornstein-Uhlenbeck semigroups on
infinite-dimensional spaces; see G\"{o}tze \cite{Goetze} for a
version of Stein's method based on multi-dimensional
Ornstein-Uhlenbeck semigroups). The \textsl{infinitesimal
generator of the Ornstein-Uhlenbeck semigroup} is noted $L$. A
square integrable random variable $F$ is in the domain of $L$
(noted ${\Dom}L$) if $F$ belongs to the domain of $\delta D$ (that
is, if $F$ is in $\mathbb{D}^{1,2}$ and $DF \in {\rm Dom}\delta$)
and, in this case,
$L F =- \delta D F .$
One can prove that $LF$ is such that
$ LF=- \sum_{n=0}^\infty n J_n(F).$
As an example, if $F=I_q(f_q)$, with $f_q \in \EuFrak{H}^{\odot
q}$, then $LF=-q F$.  Note that, for every $F\in {\Dom}L$, one has
$E(LF)=0$.
The inverse $L^{-1}$ of the operator $L$ acts on
zero-mean random variables $F\in L^2(X)$
as
$ L^{-1}F=-\sum_{n=1}^\infty \frac1n J_n(F).$
In particular, for every $q\geq 1$ and every $F=I_q(f_q)$ with
$f_q \in \EuFrak{H}^{\odot q}$, one has that $L^{-1}F=-\frac1q F$.

\smallskip

We conclude this section by
recalling
two important characterizations of the Ornstein-Uhlenbeck semigroup and its
generator.

\smallskip

\noindent i) \underline{Mehler's formula.} Let $F$ be an element
of $L^2(X)$, so that $F$ can be represented as an
application from $\real ^ \HH$ into $\real$. Then, an alternative
representation (due to Mehler) of the action of the
Ornstein-Uhlenbeck semigroup $T$ (as defined in (\ref{OrnstU})) on
$F$, is the following:
\begin{equation}
\label{MehlerF}
T_t(F) =E[F(e^{-t} a + \sqrt{1-e^{-2t}} X)] \mid _{a = X}, \quad
t\geq 0,
\end{equation}
where $a$ designs a generic element of $\real ^ \HH$. See Nualart
\cite[Section 1.4.1]{Nualart} for more details on this and other
characterizations of $T$.

\smallskip

\noindent ii) \underline{Differential characterization of $L$.}
Let $F\in L^2(X)$ have the form
$
F = f(X(h_1),...,X(h_d)),
$
where $f \in \mathscr{C}^2(\real ^d)$
has bounded first and second derivatives,
and $h_i \in \HH$, $i =
1,...,d$. Then,
\begin{eqnarray}
\label{DiffELLE}
LF& =&  \sum_{i,j=1}^d \frac{\partial^2 f}{\partial x_i \partial
x_j}(X(h_1),...,X(h_d))\langle h_i , h_j \rangle_\HH
 -\sum_{i=1}^d \frac{\partial f}{\partial
 x_i}(X(h_1),...,X(h_d))X(h_i). \quad\quad
\end{eqnarray}
See Propositions 1.4.4 and 1.4.5 in \cite{Nualart} for a proof and
some generalizations of (\ref{DiffELLE}).
\section{Stein's method and Gaussian
vectors}\label{S : Main}

We start by giving a definition of the Wasserstein distance, as
well as by introducing some useful norms over classes of
real-valued matrices.

\begin{definition}\label{DefWass}
{\rm
\begin{itemize}
\item[(i)] The {\it Wasserstein distance} between the laws of two $\real^d$-valued
random vectors $X$ and $Y$, noted $d_{\rm W}(X,Y)$, is given by
$$ d_{\rm W}(X,Y):=\sup_{g\in\mathscr{H}; \|g\|_{Lip}\leq 1} \big\vert E[g(X)]-E[g(Y)] \big\vert,$$
where $\mathscr{H}$ indicates the class of Lipschitz functions,
that is, the collection of all functions $g:\real^d\to\real$ such
that $\displaystyle{\|g\|_{Lip}:=\sup_{x\neq y}\frac{\vert
g(x)-g(y) \vert}{\| x-y \|_{\real^d}}<\infty}$ (with
$\|\cdot\|_{\real^d}$ the usual Euclidian norm on $\real^d$).
\item[(ii)] The {\it Hilbert-Schmidt inner product} and the {\it Hilbert-Schmidt norm}
on the class of $d\times d$ real matrices, denoted respectively by
$\langle \cdot, \cdot \rangle_{H.S.}$ and $\|\cdot\|_{H.S.}$, are
defined as follows: for every pair of matrices $A$ and $B$,
$ \langle A, B \rangle_{H.S.}:={\rm Tr}(A B^{T})$ and $\|A\|_{H.S.}:=\sqrt{\langle A, A \rangle_{H.S.}}.$
\item[(iii)] The {\it operator norm} of a $ d\times d$ matrix $A$ over $\real$ is given
by
$ \|A\|_{op} :=\sup_{\|x\|_{\real^d}=1}\|A
x\|_{\real^d}.$
\end{itemize}
}
\end{definition}

\begin{remark}
{\rm
\begin{itemize}
\item[1.] For every $d\geq 1$ the topology induced by $d_{\rm
W}$, on the class of all probability measures on $\real^d$, is
strictly stronger than the topology induced by weak convergence
(see e.g. Dudley \cite[Chapter 11]{Dudley book}).
\item[2.] The reason why we focus on the Wasserstein distance is nested in the statement of the
forthcoming Lemma \ref{lemma:SteinIBP2}. Indeed, according to relation (\ref{SteinEstimate}), in order to
control the second derivatives of the solution of the Stein equation (\ref{eq:SteinEquation2})
associated with $g$, one must use the fact that $g$ is Lipschitz. 
\item [3.] According to the notation introduced in Definition
\ref{DefWass}(ii), relation (\ref{DiffELLE}) can be rewritten as
\begin{equation}
\label{compactDIFF}
LF = \langle C , {\rm Hess} f(Z) \rangle_{H.S.} - \langle Z ,
\nabla f (Z) \rangle_{\real^d},
\end{equation}
where $Z =(X(h_1),...,X(h_d))$, and $C = \{C(i,j) : i,j=1,...,d\}$
is the $d \times d$ covariance matrix such that $C(i,j)
=E(X(h_i)X(h_j))=\langle h_i , h_j \rangle_\HH$.
\end{itemize}
}
\end{remark}

\smallskip
Given a $d\times d$ positive definite symmetric matrix $C$, we use
the notation $\mathscr{N}_d(0,C)$ to indicate the law of a
$d$-dimensional Gaussian vector with zero mean and covariance $C$.
The following result, which is basically known (see e.g.
\cite{ChatterjeeMeckes} or \cite{ReinertRollin}), is the
$d$-dimensional counterpart of Stein's Lemma
\ref{lemma:SteinIBP1}. In what follows, we provide a new proof
which is almost exclusively based on the use of Malliavin
operators.

\begin{lemma}
\label{lemma:SteinIBP2} Fix an integer $d\geq 2$ and let
$C=\{C(i,j) : i,j=1,...,d\}$ be a $d \times d$ positive definite
symmetric real matrix.
\begin{itemize}
\item[(i)] Let $Y$ be a
random variable with values in $\real^d$. Then
$Y\sim\mathscr{N}_d(0,C)$ if and only if, for every twice
differentiable function $f:\real^d\to\real$ such that $ E\vert
\langle C , {\rm Hess} f(Y) \rangle_{H.S.}\vert+E\vert \langle Y, \nabla f(Y)
\rangle_{\real^d}\vert <\infty$, it holds that
\begin{equation}
\label{SteinCharD}
 E[
\langle Y, \nabla f(Y)
\rangle_{\real^d}
-\langle C,
{\rm Hess} f(Y) \rangle_{H.S.}
]=0.
\end{equation}
\item[(ii)] Consider a Gaussian random vector $Z\sim \mathscr{N}_d(0,C)$.
Let $g:\real^d\to\real$ belong to $\mathscr{C}^2(\real^d)$
with first and second bounded derivatives.
Then, the function $U_0(g)$ defined by
\begin{equation*}
U_0g(x):=\int_0^1 \frac{1}{2 t}
E[g(\sqrt{t}x+\sqrt{1-t}Z)-g(Z)] dt
\end{equation*}
is a solution to the following differential equation (with unknown function $f$):
\begin{equation}
\label{eq:SteinEquation2}
g(x)-E[g(Z)]=
\langle x, \nabla f(x) \rangle_{\real^d}
-\langle C , {\rm Hess} f(x)
\rangle_{H.S.}
, \quad
x\in\real^d.
\end{equation}
Moreover, one has that
\begin{equation}
\label{SteinEstimate}
\sup_{x\in \real^d} \| {\rm Hess}\,U_0g(x) \|_{H.S.}\leq \|
C^{-1}\|_{op} \,\, \| C\|_{op}^{1/2} \,\, \|g\|_{Lip}.
\end{equation}
\end{itemize}
\end{lemma}

\begin{remark}
{\rm
\begin{itemize}
\item[1.] If $C = \sigma^2 \textbf{I}_d$ for some $\sigma >0$ (that is,
if $Z$ is composed of i.i.d. centered Gaussian random variables
with common variance equal to $\sigma^2$), then
$$
\| C^{-1}\|_{op} \,\, \| C\|_{op}^{1/2} = \| \sigma^{-2}
\textbf{I}_d\|_{op} \,\, \| \sigma^2 \textbf{I}_d \|_{op}^{1/2} =
\sigma^{-1}.
$$
\item[2.] Unlike formulae (\ref{urStein}) and
(\ref{eq:SteinEquation1}) (associated with one-dimensional
Gaussian approximations) the relation (\ref{SteinCharD}) and the
Stein equation (\ref{eq:SteinEquation2}) involve second-order
differential operators. A discussion of this fact is detailed e.g.
in \cite[Theorem 4]{ChatterjeeMeckes}.
\end{itemize}
}
\end{remark}

\noindent
{\it Proof of Lemma \ref{lemma:SteinIBP2}.} We start by proving
Point (ii). First observe that, without loss of generality, we can
suppose that $Z=(Z_1,...,Z_d):=(X(h_1),...X(h_d))$, where $X$ is an
isonormal Gaussian process over
$\HH=\real^d$,
the
kernels $h_i $ belong to $\HH$ ($i=1,...,d$), and $\langle h_i ,
h_j \rangle_\HH  = E(X(h_i)X(h_j))= E(Z_i Z_j) = C(i,j)$. By using
the change of variable $2u =-\log t$, one can rewrite $U_0 g(x)$
as follows
$$
U_0 g(x) = \int_0^\infty \{E[g(e^{-u}
x+\sqrt{1-e^{-2u}}Z)]-E[g(Z)]\} du.
$$
Now define $\widetilde{g}(Z) := g(Z)-E[g(Z)]$, and observe that
$\widetilde{g}(Z)$ is by assumption a centered element of ${\rm
L}^2(X)$. For $q\geq 0$, denote by $J_q(\widetilde{g}(Z))$ the
projection of $\widetilde{g}(Z)$ on the $q$th Wiener chaos, so
that $J_0(\widetilde{g}(Z))=0$. According to Mehler's formula
(\ref{MehlerF}),
$$E[g(e^{-u} x+\sqrt{1-e^{-2u}}Z)]|_{x = Z}-E[g(Z)]  =
E[\widetilde{g}(e^{-u} x+\sqrt{1-e^{-2u}}Z)]|_{x = Z} = T_u
\widetilde{g}(Z),$$ where $x$ denotes a generic element of $\real
^d$. In view of (\ref{OrnstU}), it follows that
$$
U_0 g(Z) = \int_0^\infty T_u \widetilde{g}(Z) du = \int_0^\infty
\sum_{q\geq 1} e^{-qu} J_q(\widetilde{g}(Z)) du = \sum_{q\geq 1}
\frac1q J_q(\widetilde{g}(Z)) = -L^{-1}\widetilde{g}(Z).
$$
Since $g$ belongs to $\mathscr{C}^2(\real^d)$ with bounded first and second derivatives,
it is easily seen that the same holds for $U_0 g$.
By exploiting the differential
representation (\ref{compactDIFF}),
one deduces that
$$
\langle Z ,
\nabla U_0 g (Z) \rangle_{\real^d}
-
\langle C , {\rm Hess} U_0 g(Z)  \rangle_{H.S.}
= -L U_0 g(Z) =
LL^{-1}\widetilde{g}(Z) = g(Z) - E[g(Z)].
$$
Since the matrix $C$ is positive definite, we infer that the
support of the law of $Z$ coincides with $\real ^d$, and therefore
(e.g. by a continuity argument) we obtain that
$$
\langle x ,
\nabla U_0 g (x) \rangle_{\real^d}
-
\langle C , {\rm Hess}\, U_0 g(x)  \rangle_{H.S.}
= g(x) - E[g(Z)],
$$
for every $x\in\real ^d$. This yields that the function $U_0g$
solves the Stein's equation (\ref{eq:SteinEquation2}).

To prove the estimate (\ref{SteinEstimate}), we first
recall that there exists a unique non-singular \textsl{symmetric}
matrix $A$ such that $A^2 = C$, and that one has that $A^{-1} Z
\sim \mathscr{N}_d(0,\mathbf{I}_d)$. Now write $U_0g(x) =
h(A^{-1}x)$, where
$$
 h(x)= \int_0^1 \frac{1}{2 t}
E[g_A(\sqrt{t}x+\sqrt{1-t}A^{-1}Z)-g_A(A^{-1}Z)] dt,
$$
and $g_A(x)=g(Ax)$. Note that, since $A^{-1} Z \sim
\mathscr{N}_d(0,\mathbf{I}_d)$, the function $h$ solves the Stein's
equation
$
\langle x , \nabla h (x) \rangle_{\real^d}-\Delta h(x)  = g_A(x)
- E[g_A(Y)],
$
where $Y\sim \mathscr{N}_d(0,\mathbf{I}_d)$. We can now use the same
arguments as in the proof of Lemma 3 in \cite{ChatterjeeMeckes}
to deduce that
\begin{equation}
\label{MeckessMeckess}
\sup_{x\in \real^d} \| {\rm Hess}\,h(x) \|_{H.S.}\leq
 \|g_A\|_{Lip}\leq \|A\|_{op}\|g\|_{Lip}.
\end{equation}
On the other hand, by noting $h_{A^{-1}}(x) = h(A^{-1}x)$, one
obtains by standard computations (recall that $A$ is symmetric)
that
$
{\rm Hess}\,U_0g(x)={\rm Hess}\,h_{A^{-1}}(x) = A^{-1}{\rm
Hess}\,h(A^{-1}x) A^{-1},
$
yielding
\begin{eqnarray}
\sup_{x\in \real^d} \| {\rm Hess}\,U_0g(x) \|_{H.S.}&=& \sup_{x\in
\real^d} \|A^{-1}  {\rm Hess}\,h(A^{-1}x) A^{-1} \|_{H.S.} \notag \\
&=& \sup_{x\in
\real^d} \|A^{-1}  {\rm Hess}\,h(x) A^{-1} \|_{H.S.} \notag \\
& \leq & \|A^{-1}\|_{op}^2 \sup_{x\in \real^d} \| {\rm
Hess}\,h(x) \|_{H.S.} \label{Sing1}\\
& \leq &\|A^{-1}\|_{op}^2 \,\, \| A \|_{op} \,\, \|g\|_{Lip} \label{Sing2} \\
& \leq &\| C^{-1}\|_{op} \,\, \| C\|_{op}^{1/2} \,\, \|g\|_{Lip}.
\label{Sing3}
\end{eqnarray}
The chain of inequalities appearing in formulae
(\ref{Sing1})--(\ref{Sing3}) are mainly a consequence of the usual
properties of the Hilbert-Schmidt and operator norms. Indeed, to
prove inequality (\ref{Sing1}) we used the relations
\begin{eqnarray*}  \|A^{-1}  {\rm Hess}\,h(x) A^{-1}
\|_{H.S.}& \leq &
\|A^{-1}\|_{op} \,\, \| {\rm Hess}\,h(x) A^{-1}\|_{H.S.} \\
&\leq & \|A^{-1}\|_{op}\,\, \| {\rm Hess}\,h(x)
\|_{H.S.}\,\,\|A^{-1}\|_{op}\,\,;
\end{eqnarray*}
relation (\ref{Sing2}) is a consequence of (\ref{MeckessMeckess});
finally, to show the inequality (\ref{Sing3}), one uses the fact
that
\begin{eqnarray*}
\|A^{-1}\|_{op}  \leq  \sqrt{\|A^{-1} A^{-1}\|_{op}}=
\sqrt{\|C^{-1}\|_{op}} \quad\mbox{and}\quad
\|A\|_{op}  \leq  \sqrt{\|A A\|_{op}}= \sqrt{\|C\|_{op}} \,\,.
\end{eqnarray*}
\noindent We are now left with the proof of Point (i) in the
statement. The fact that a vector $Y\sim\mathscr{N}_d(0,C)$
necessarily verifies (\ref{SteinCharD}) can be proved by standard
integration by parts. On the other hand, suppose that $Y$ verifies
(\ref{SteinCharD}). Then, according to Point (ii), for every $g\in
\mathscr{C}^2(\real^d)$
with bounded first and second derivatives, 
$$
E(g(Y)) - E(g(Z)) =E(\langle Y, \nabla U_0g(Y) \rangle_{\real^d} - \langle C , {\rm Hess}\, U_0g(Y)
\rangle_{H.S.}) = 0,
$$
where $Z\sim \mathscr{N}_d(0,C)$. Since the collection of all such
functions $g$ generates the Borel $\sigma$-field on $\real^d$,
this implies that $Y \stackrel{\rm Law}{=} Z$, thus yielding the
desired conclusion.
\qed \\

The following statement is the main result of this paper. Its
proof makes a crucial use of the integration by parts formula
(\ref{eq:MalliavinIBP}) discussed in Section
\ref{section:Malliavin}.

\begin{theorem}
\label{theo:majDist} Fix $d\geq 2$ and let $C = \{C(i,j) :
i,j=1,...,d\}$ be a $d\times d$ positive definite matrix. Suppose
that $Z\sim\mathscr{N}_d(0,C)$ and that $F=(F_1,\ldots,F_d)$ is a
$\real^d$-valued random vector such that $E[F_i]=0$ and $F_i \in
\mathbb{D}^{1,2}$ for every $i=1,\ldots,d$. Then,
\begin{eqnarray}
 d_{\rm W}(F,Z)&\leq &\| C^{-1}\|_{op} \,\,
\| C\|_{op}^{1/2} \sqrt{E\| C - \Phi(D F) \|_{H.S}^2}\label{WassONE}
\\
&=& \| C^{-1}\|_{op} \,\, \| C\|_{op}^{1/2} \sqrt{\sum_{i,j=1}^d
E[( C(i,j) - \langle DF_i, -DL^{-1}F_j \rangle_{\EuFrak H})^2] },
\label{WassTWO}
\end{eqnarray}
where we write $\Phi(DF)$ to indicate the matrix
$\Phi(DF):=\{\langle DF_i, -DL^{-1}F_j \rangle_{\EuFrak H}:1\leq
i,j \leq d \}.$
\end{theorem}

\begin{proof}
We start by proving that, for every $g\in\mathscr{C}^2(\real^d)$
with bounded first and second derivatives,
\begin{equation*}
|E[g(F)] - E[g(Z)]| \leq \| C^{-1}\|_{op} \,\, \|
C\|_{op}^{1/2}\,\, \|g\|_{Lip} \sqrt{E\| C - \Phi(D F)
\|_{H.S}^2}.
\end{equation*}
To prove such a claim, observe that, according to Point (ii) in
Lemma \ref{lemma:SteinIBP2}, $E[g(F)] - E[g(Z)] =
E[
 \langle F , \nabla U_0 g (F)
\rangle_{\real^d}
-
\langle C ,
{\rm Hess} U_0 g(F)  \rangle_{H.S.}
]$. Moreover,
\begin{eqnarray*}
&&\displaystyle{ \big\vert E[\langle C , {\rm Hess} U_0 g(F)
\rangle_{H.S.} - \langle F , \nabla U_0 g (F)
\rangle_{\real^d}]  \big\vert}\\
&=&\displaystyle{\left\vert
E\left[ \sum_{i,j=1}^d C(i,j)\partial^2_{ij} U_0g(F) - \sum_{i=1}^d F_i \partial_i U_0g(F)\right] \right\vert}\\
&=&\displaystyle{\left\vert\sum_{i,j=1}^d E\left[
C(i,j)\partial^2_{ij} U_0g(F)\right] -
\sum_{i=1}^d E\left[\big(LL^{-1}F_i\big) \partial_i U_0g(F)\right]  \right\vert} \,\,\, (\textrm{since } E(F_i)=0) \\
&=&\displaystyle{\left\vert \sum_{i,j=1}^d E\left[C(i,j)
\partial^2_{ij} U_0 g(F)\right] +
\sum_{i=1}^d E\left[\delta(DL^{-1}F_i) \partial_i U_0g(F)\right]   \right\vert\,\,\, (\textrm{since } \delta D=-L})\\
&=&\displaystyle{ \left\vert\sum_{i,j=1}^d  E\left[C(i,j)
\partial^2_{ij} U_0g(F)\right]\!\!- \!\!\sum_{i=1}^d E\left[\langle
D(\partial_i U_0g(F)),-DL^{-1}F_i\rangle_\EuFrak{H}\right]  \right\vert \,\,\,\, (\textrm{by } (\ref{eq:MalliavinIBP}) })\\
&=&\displaystyle{ \left\vert\sum_{i,j=1}^d  E\left[C(i,j)
\partial^2_{ij} U_0g(F)\right] - \sum_{i,j=1}^d E\left[\partial^2_{ji}
U_0g(F)  \langle  DF_j,-DL^{-1}F_i \rangle_\EuFrak{H}\right]  \right\vert}\,\,\,\, (\textrm{by }(\ref{eq:ChainRule})) \\
&=&\displaystyle{ \left\vert\sum_{i,j=1}^d E\left[\partial^2_{ij}
U_0g(F)  \big( C(i,j)-
\langle  DF_i,-DL^{-1}F_j \rangle_\EuFrak{H} \big)\right] \right\vert }\\
&=&\displaystyle{ \big\vert E\langle {\rm Hess}\, U_0g(F),C- \Phi(DF)\rangle_{H.S.} \big\vert}\\
&\leq& \displaystyle{\sqrt{E\|{\rm Hess} \,U_0g(F)
\|_{H.S}^2} \sqrt{E\| C - \Phi(D F) \|_{H.S}^2}\,\,\, (\textrm{by the Cauchy-Schwarz inequality}})\\
&\leq& \| C^{-1}\|_{op} \,\, \| C\|_{op}^{1/2}\,\, \|g\|_{Lip}
\sqrt{E\| C - \Phi(D F) \|_{H.S}^2} \,\,\,(\textrm{by
(\ref{SteinEstimate})}).
\end{eqnarray*}
To prove the Wasserstein estimate (\ref{WassONE}), it is
sufficient to observe that, for every globally Lipschitz function
$g$ such that $\|g\|_{Lip}\leq 1$, there exists a
family $\{g_\e
: \e> 0\}$ such that:
\begin{itemize}
\item[(i)] for each $\e>0$, the first and second derivatives of $g_\e$
are bounded;
\item[(ii)] for each $\e>0$, one has that $\|g_\e\|_{Lip} \leq \|g\|_{Lip}$;
\item[(iii)] as $\e\rightarrow 0$, $\|g_\e - g\|_\infty \downarrow 0$.
\end{itemize}
For instance, we can choose $g_\e(x)=E\big[g(x+\sqrt{\e}N)\big]$ with $N\sim\mathscr{N}_d(0,{\bf I}_d)$.
\end{proof}

Observe that Theorem \ref{theo:majDist} generalizes relation
(\ref{GioIvan}) (that was proved in \cite[Theorem
3.1]{NourdinPeccati1}). We now aim at applying Theorem
\ref{theo:majDist} to vectors of multiple stochastic integrals.

\begin{corollary}
\label{prop:Chaos} Fix $d\geq 2$ and $1\leq q_1\leq\ldots\leq
q_d$. Consider a vector
$F:=(F_1,\ldots,F_d)=(I_{q_1}(f_1),\ldots,I_{q_d}(f_d))$ with
$f_{i}\in \EuFrak{H}^{\odot q_i}$ for any $i=1\ldots,d$. Let
$Z\sim\mathscr{N}_d(0,C)$, with $C$ positive definite. Then,
\begin{equation}
\label{star}
d_{\rm W} (F,Z) \leq \| C^{-1}\|_{op} \,\, \| C\|_{op}^{1/2}
 \sqrt{ \! \sum_{1\le i, j\le d} E\left[\left( C(i,j)\!-\!
\frac{1}{q_j} \langle DF_i,DF_j \rangle_{\EuFrak{H}}\right)^2
\right]}.
\end{equation}
\end{corollary}

\begin{proof}
We have $-L^{-1}F_j=\frac{1}{q_j}\,F_j$ so that
the desired conclusion follows from (\ref{WassTWO}).
\end{proof}

When one applies Corollary \ref{prop:Chaos} in concrete situations
(see e.g. Section \ref{section:applications} below), one can use
the following result in order to evaluate the 
right-hand side
of (\ref{star}).
\begin{lemma}\label{lm-control}
Let $F=I_p(f)$ and $G=I_q(g)$, with $f\in\HH^{\odot p}$ and
$g\in\HH^{\odot q}$ ($p,q\geq 1$). Let $a$ be a real constant. If
$p=q$, one has the estimate:
\begin{eqnarray*}
&&E\left[\left(a-\frac1p\left\langle DF,DG\right\rangle_\HH\right)^2\right] \leq (a-p!\langle f,g\rangle_{\HH^{\otimes p}})^2\\
&&\hskip2cm+
\frac{p^2}{2}\sum_{r=1}^{p-1}(r-1)!^2\binom{p-1}{r-1}^4(2p-2r)!
\big( \|f\otimes_{p-r}f\|^2_{\HH^{\otimes
2r}}+\|g\otimes_{p-r}g\|^2_{\HH^{\otimes 2r}}\big).
\end{eqnarray*}
On the other hand, if $p< q$, one has that
\begin{eqnarray*}
&&E\left[\left(a-\frac1q\left\langle DF,DG\right\rangle_\HH
\right)^2\right]
\leq a^2+p!^2\binom{q-1}{p-1}^2(q-p)!\|f\|^2_{\HH^{\otimes p}}\|g\otimes_{q-p}g\|_{\HH^{\otimes 2p}}\\
&&+
\frac{p^2}{2}\sum_{r=1}^{p-1}(r-1)!^2\binom{p-1}{r-1}^2\binom{q-1}{r-1}^2(p+q-2r)!\big(
\|f\otimes_{p-r}f\|^2_{\HH^{\otimes
2r}}+\|g\otimes_{q-r}g\|^2_{\HH^{\otimes 2r}}\big).
\end{eqnarray*}
\end{lemma}
\begin{remark}\label{r}
{\rm
\begin{enumerate}
\item Recall that $E\big(I_p(f)I_q(g)\big)=\left\{\begin{array}{lll}
p!\langle f,g\rangle_{\HH^{\otimes p}}&\quad\mbox{if $p=q$},\\
0&\quad\mbox{otherwise}.
\end{array}\right.$
\item In order to estimate the right-hand side of (\ref{star}), we see that it
is sufficient to asses the quantity $\|f_i\otimes_r
f_i\|_{\HH^{\otimes2(q_i-r)}}$ for any $i\in\{1,\ldots,d\}$ and
$r\in\{1,\ldots,q_i-1\}$ on the one hand, and $\langle
f_i,f_j\rangle_{\HH^{\otimes q_i}}$ for any $1\leq i,j\leq d$ such
that $q_i=q_j$ on the other hand.
\end{enumerate}
}
\end{remark}
\noindent {\it Proof of Lemma \ref{lm-control}} (see also \cite[Lemma 2]{NualartOrtizLatorre}). Without loss of
generality, we can assume that $\HH=L^{2}(A,\mathscr{A},\mu)$,
where $(A,\mathscr{A})$ is a measurable space, and $\mu$ is a
$\sigma$-finite and non-atomic measure. Thus, we can write
\begin{eqnarray*}
\langle DF,DG\rangle_\HH &=&p\,q\left\langle I_{p-1}(f),I_{q-1}(g)\right\rangle_\HH
=p\,q\int_A I_{p-1}\big(f(\cdot,t)\big)I_{q-1}\big(g(\cdot,t)\big)\mu(dt)\\
&=&p\,q\int_A \sum_{r=0}^{p\wedge q-1} r!\binom{p-1}{r} \binom{q-1}{r} I_{p+q-2-2r}\big(f(\cdot,t)\widetilde{\otimes}_r g(\cdot,t)\big)\mu(dt)\\
&=&p\,q\sum_{r=0}^{p\wedge q-1} r!\binom{p-1}{r}\binom{q-1}{r} I_{p+q-2-2r}(f\widetilde{\otimes}_{r+1}g)\\
&=&p\,q \sum_{r=1}^{p\wedge q} (r-1)!\binom{p-1}{r-1}\binom{q-1}{r-1} I_{p+q-2r}(f\widetilde{\otimes}_r g).
\end{eqnarray*}
It follows that
\begin{eqnarray}
&&E\left[\left(a-\frac1q\left\langle DF,DG\right\rangle_\HH\right)^2\right] \label{Murray}\\
&=&\left\lbrace
\begin{array}{l}
a^2+p^2\sum_{r=1}^{p}(r-1)!^2
\binom{p-1}{r-1}^2\binom{q-1}{r-1}^2 (p+q-2r)!
\|f\widetilde{\otimes}_r g\|^2_{\HH^{\otimes (p+q-2r)}} \textrm{ if } p< q,\\\\
(a-p!\langle f, g\rangle_{\EuFrak{H}^{\otimes
p}})^2+p^2\sum_{r=1}^{p-1}(r-1)!^2 \binom{p-1}{r-1}^4 (2p-2r)!
\|f\widetilde{\otimes}_r g\|^2_{\HH^{\otimes (2p-2r)}} \textrm{ if
} p=q.
\end{array}\notag
\right.
\end{eqnarray}
If $r<p\leq q$ then
\begin{eqnarray*}
\|f\widetilde{\otimes}_r g\|^2_{\HH^{\otimes (p+q-2r)}}
&\leq& \|f\otimes_r g\|^2_{\HH^{\otimes (p+q-2r)}}
=\langle f\otimes_{p-r} f, g\otimes_{q-r}g\rangle_{\HH^{\otimes 2r}}\\
&\leq&
\|f\otimes_{p-r}f\|_{\HH^{\otimes 2r}}\|g\otimes_{q-r}g\|_{\HH^{\otimes 2r}}\\
&\leq&\frac12\left(
\|f\otimes_{p-r}f\|_{\HH^{\otimes 2r}}^2+\|g\otimes_{q-r}g\|_{\HH^{\otimes 2r}}^2
\right).
\end{eqnarray*}
If $r=p<q$, then
$$
\|f\widetilde{\otimes}_p\, g\|^2_{\HH^{\otimes (q-p)}} \leq
\|f\otimes_p \,g\|^2_{\HH^{\otimes (q-p)}} \leq
\|f\|^2_{\HH^{\otimes p}}\|g\otimes_{q-p}g\|_{\HH^{\otimes 2p}}.
$$
If $r=p=q$, then $ f\widetilde{\otimes}_p g =\langle
f,g\rangle_{\HH^{\otimes p}}. $
By plugging these last expressions into (\ref{Murray}), we deduce
immediately the desired conclusion.
\qed\\

Let us now recall the following result, which is a collection of
some of the findings contained in the papers by Peccati and Tudor
\cite{PeccatiTudor} and Nualart and Ortiz-Latorre
\cite{NualartOrtizLatorre}.

\begin{theorem}[See \cite{NualartOrtizLatorre,PeccatiTudor}]
\label{theo:NOL-PT} Fix $d\geq2$ and let
$C=\{C(i,j):i,j=1,...,d\}$ be a $d\times d$ positive definite
matrix. Fix integers $1\leq q_1\leq\ldots\leq q_d$. For any $n\geq
1$ and $i=1,\ldots,d$, let $f_{i}^{(n)}$ belong to $
\EuFrak{H}^{\odot q_i}$. Assume that
$$F^{(n)}=(F^{(n)}_1,\ldots,F^{(n)}_d):=(I_{q_1}(f_{1}^{(n)}),
\ldots,I_{q_d}(f_{d}^{(n)}))\quad n\geq 1,$$ is such that
\begin{equation}
\label{eq:asympcov}
\lim_{n\to\infty}
E[F_i^{(n)}F_j^{(n)}]=C(i,j),\quad 1\leq i,j\leq d.
\end{equation}
Then, as $n\to\infty$, the following four assertions are
equivalent:
\begin{itemize}
\item[(i)] For every $1\leq i\leq d$, $F_i^{(n)}$ converges in distribution
to a centered Gaussian random variable with variance $C(i,i)$.
\item[(ii)] For every $1\leq i\leq d$, $E\left[(F_i^{(n)})^4\right]\rightarrow3C(i,i)^2$.
\item[(iii)] For every $1\leq i\leq d$ and every $1\leq r \leq q_i-1$,
$\|f_{i}^{(n)} \otimes_r f_{i}^{(n)}\|_{\EuFrak{H}^{\otimes
2(q_i-r)}} \rightarrow 0$.
\item[(iv)] The vector $F^{(n)}$ converges in distribution to a $d$-dimensional Gaussian
vector $\mathscr{N}_d(0,C)$.
\end{itemize}
Moreover, if $C(i,j)= \delta_{ij}$, where $\delta_{ij}$ is the
Kronecker symbol, then either one of conditions (i)--(iv) above is
equivalent to the following:
\begin{itemize}
\item[(v)] For every $1\leq i\leq d$, $\|DF_i^{(n)}\|_{\EuFrak{H}}^2 \overset{L^2}{\longrightarrow} q_i$.
\end{itemize}
\end{theorem}

We conclude this section by pointing out the remarkable fact that,
for vectors of multiple Wiener-It\^{o} integrals of arbitrary
length, \textsl{the Wasserstein distance metrizes the weak
convergence towards a Gaussian vector with positive definite
covariance}. Note that the next statement also contains a
generalization of the multidimensional results proved in
\cite{NualartOrtizLatorre} to the case of an arbitrary covariance.

\begin{prop}
Fix $d\geq 2$, let $C$ be a positive definite $d\times d$
symmetric matrix, and let $1\leq q_1\leq\ldots\leq q_d$. Consider
vectors
$$F^{(n)}:=(F_1^{(n)},\ldots,F_d^{(n)})=
(I_{q_1}(f_1^{(n)}),\ldots,I_{q_d}(f^{(n)}_d)), \quad n\geq 1,$$
with $f^{(n)}_{i}\in \EuFrak{H}^{\odot q_i}$ for every
$i=1\ldots,d$. Assume moreover that $F^{(n)}$ satisfies condition
(\ref{eq:asympcov}). Then, as $n\to\infty$,  the following three
conditions are equivalent:
\begin{itemize}
\item[(a)] $d_{\rm W}(F^{(n)},Z)\rightarrow 0$.
\item[(b)] For every $1\leq i\leq d$, $q_i^{-1}\|DF_i^{(n)}\|_{\EuFrak{H}}^2
\overset{\rm L^2}{\longrightarrow}C(i,i)$ and, for every  $1\leq
i\neq j \leq d$, $\langle DF^{(n)}_i,-DL^{-1}F^{(n)}_j \rangle_{\HH} =
q_j^{-1}\langle DF^{(n)}_i,DF^{(n)}_j \rangle_{\HH} \overset{\rm
L^2}{\longrightarrow}C(i,j)$.
\item[(c)] $F^{(n)}$
converges in distribution to $Z\sim\mathscr{N}_d(0,C)$.
\end{itemize}
\end{prop}

\begin{proof}
Since convergence in the Wasserstein distance implies convergence
in distribution, the implication (a) $\rightarrow$ (c) is trivial.
The implication (b) $\rightarrow$ (a) is a consequence of relation
(\ref{star}). Now assume that (c) is verified, that is, $F^{(n)}$
converges in law to $Z\sim\mathscr{N}_d(0,C)$ as $n$ goes to
infinity. By Theorem \ref{theo:NOL-PT} we have that, for any
$i\in\{1,\ldots,d\}$ and $r\in\{1,\ldots,q_i-1\}$,
$$
\|f_i^{(n)}\otimes_r f_i^{(n)}\|_{\HH^{\otimes2(q_i-r)}} \underset{n\to\infty}{\longrightarrow} 0.
$$
By combining Corollary \ref{prop:Chaos} with Lemma \ref{lm-control}
(see also Remark \ref{r}(2)), one therefore easily deduces that,
since (\ref{eq:asympcov}) is in order, condition (b) must
necessarily be satisfied.
\end{proof}

\section{Applications} \label{section:applications}

\subsection{Convergence of marginal distributions in the functional
Breuer-Major CLT} \label{SS : BreuMaj}
 In this section, we
use our main results in order to derive an explicit bound for the
celebrated \textsl{Breuer-Major CLT} for fractional Brownian
motion (fBm). We recall that a fBm $B=\{B_t:t\geq 0\}$, with Hurst
index $H\in(0,1)$, is a centered Gaussian process, started from zero and with covariance function
$E(B_{s}B_{t})=R(s,t)$, where
$$
R(s,t)=\frac{1}{2}\left( t^{2H}+s^{2H}-|t-s|^{2H}\right);\quad
s,t\geq 0.
$$
For any choice of the Hurst parameter
$H\in(0,1)$, the Gaussian space generated by $B$ can be identified
with an isonormal Gaussian process of the type
$X=\{X(h):h\in\HH\}$, where the real and separable Hilbert space
$\EuFrak H$ is defined as follows: (i) denote by $\mathscr{E}$ the
set of all $\mathbb{R}$-valued step functions on $[0,\infty)$, (ii)
define $\EuFrak H$ as the Hilbert space obtained by closing
$\mathscr{E}$ with respect to the scalar product
$
\left\langle
{\mathbf{1}}_{[0,t]},{\mathbf{1}}_{[0,s]}\right\rangle _{\EuFrak
H}=R(t,s).
$
In particular, with such a notation, one has that
$B_t=X(\mathbf{1}_{[0,t]})$.
The reader is referred e.g. to \cite{Nualart} for more details on fBm, including crucial connections
with fractional operators. We also define $\rho(\cdot)$ to be the
covariance function associated with the stationary process
$x\mapsto B_{x+1}-B_x$ ($x\in\mathbb{R}$), that is
$$
\rho(x):=\frac12\big(|x+1|^{2H}+|x-1|^{2H}-2|x|^{2H}\big)\underset{|x|\to\infty}{\sim}
H|2H-1|\,|x|^{2H-2}.
$$
Now fix an integer $q\geq 2$, assume that $H<1-\frac1{2q}$ and set
$$
S_n(t)=\frac{1}{\sigma\sqrt{n}}\sum_{k=0}^{\lfloor nt\rfloor -1}H_q(B_{k+1}-B_k),\quad
t\geq 0,
$$
where $H_q$ is the $q$th Hermite polynomial defined as
\begin{equation}\label{her-pol}
H_q(x)=(-1)^q e^{x^2/2} \frac{d^q}{dx^q}e^{-x^2/2},\quad x\in\mathbb{R},
\end{equation}
and where
$
\sigma=\sqrt{q!\sum_{r\in\mathbb{Z}}\rho^2(r)}.
$
According e.g. to the main results in \cite{BM} or \cite{GS}, one
has the following CLT:
$$
\{S_n(t),\,t\geq 0\} \,\,\,\xrightarrow[n\to\infty]{{\rm
f.d.d.}}\,\,\, \mbox{standard Brownian motion},
$$
where `f.d.d.' indicates convergence in the sense of
finite-dimensional distributions. To our knowledge, the following
statement contains the first multidimensional bound for the
Wasserstein distance ever proved for $\{S_n(t),\,t\geq 0\}$.
\begin{thm}\label{BM-rev}
For any fixed $d\geq 1$ and $0=t_0<t_1<\ldots<t_d$, there exists a
constant $c$, (depending only on $d$, $H$ and $(t_0,
t_1,\ldots,t_d)$, \underline{and not} on $n$) such that, for every
$n\geq 1$:
\begin{eqnarray*}
d_{\rm W}\left(
\left(\frac{S_n(t_i)-S_n(t_{i-1})}{\sqrt{t_i-t_{i-1}}}\right)_{1\leq
i \leq d};\mathscr{N}_d(0,{\bf I}_d) \right)&\leq& c\times
\left\{\begin{array}{lll}
n^{-\frac12}&\,\,\mbox{if $H\in (0,\frac12]$}\\
\\
n^{H-1}&\,\,\mbox{if $H\in (\frac12,\frac{2q-3}{2q-2}]$}\\
\\
n^{qH-q+\frac12}&\,\,\mbox{if $H\in (\frac{2q-3}{2q-2},\frac{2q-1}{2q})$}\\
\end{array}\right..
\end{eqnarray*}
\end{thm}
\begin{proof}
Fix $d\geq 1$ and $t_0=0<t_1<\ldots<t_d$. In the sequel, $c$ will denote a constant independent of $n$, which can differ from one line to another.

First, observe that
$$\frac{S_n(t_i)-S_n(t_{i-1})}{\sqrt{t_i-t_{i-1}}}=I_q(f_i^{(n)})\quad\mbox{with}\quad
f_i^{(n)}=\frac{1}{\sigma\sqrt{n}
\sqrt{t_i-t_{i-1}}
}
\sum_{k=\lfloor nt_{i-1}\rfloor}^{\lfloor nt_i\rfloor -1}
{\bf 1}_{[k,k+1]}^{\otimes q}.
$$

In \cite{NourdinPeccati1},
proof of Theorem 4.1, it is shown that, for any $i\in\{1,\ldots,d\}$ and $r\in\{1,\ldots,q_i-1\}$:
\begin{equation}\label{bound1}
\|f_i^{(n)}\otimes_r f_i^{(n)}\|_{\HH^{\otimes 2(q_i-r)}}\leq c\times
\left\{\begin{array}{lll}
n^{-\frac12}&\,\,\mbox{if $H\in (0,\frac12]$}\\
\\
n^{H-1}&\,\,\mbox{if $H\in (\frac12,\frac{2q-3}{2q-2}]$}\\
\\
n^{qH-q+\frac12}&\,\,\mbox{if $H\in (\frac{2q-3}{2q-2},\frac{2q-1}{2q})$}\\
\end{array}\right..
\end{equation}
Moreover, when $1\leq i<j\leq d$, we have:
\begin{eqnarray}\notag
&&\big|\langle f_i^{(n)},f_j^{(n)}\rangle_{\HH^{\otimes q}}\big|\\
&=&\left|
\frac{1}{\sigma^2\,n
\sqrt{t_i-t_{i-1}}\sqrt{t_j-t_{j-1}}
}
\sum_{k=\lfloor nt_{i-1}\rfloor}^{\lfloor nt_i\rfloor -1}
\sum_{l=\lfloor nt_{j-1}\rfloor}^{\lfloor nt_j\rfloor -1}
\rho^q(l-k)\right|\notag\\
&=& \frac{c}{n}\left|
\sum_{|r|=\lfloor nt_{j-1}\rfloor-\lfloor nt_{i}\rfloor+1}
^{\lfloor nt_{j}\rfloor-\lfloor nt_{i-1}\rfloor-1}
\big[ ( \lfloor nt_{j}\rfloor -1 -r)\wedge( \lfloor nt_{i}\rfloor -1)
-( \lfloor nt_{j-1}\rfloor  -r)\vee( \lfloor nt_{i-1}\rfloor )
\big]
\rho^q(r)\right|\notag\\
&\leq& c\,\,\frac{\lfloor nt_{i}\rfloor -\lfloor nt_{i-1}\rfloor -1}{n}
\sum_{|r|\geq\lfloor nt_{j-1}\rfloor-\lfloor nt_{i}\rfloor+1}
\big|\rho(r)\big|^q
=O\big(n^{2qH-2q+1}\big),\quad\mbox{as $n\to\infty$},\label{bound2}
\end{eqnarray}
the last equality coming from
$$
\sum_{|r|\geq N}\big|\rho(r)\big|^q=O(\sum_{|r|\geq N}|r|^{2qH-2q})=O(N^{2qH-2q+1}),\quad\mbox{as $N\to\infty$}.
$$

Finally, by combining (\ref{bound1}), (\ref{bound2}), Corollary \ref{prop:Chaos} and Lemma \ref{lm-control}, we obtain
the desired conclusion.
\end{proof}

\subsection{Vector-valued functionals of finite Gaussian sequences} \label{SS : Chatterjee}

Let $Y=(Y_1,...,Y_n)\sim \mathscr{N}_n(0, {\bf I}_n)$, and let
$f:\real^n \rightarrow \real$ be an absolutely continuous function
such that $f$ and its partial derivatives have subexponential
growth at infinity. The following result has been proved by
Chatterjee in \cite{Chatterjee_ptrf}, in the context of limit
theorems for linear statistics of eigenvalues of random matrices.
We use the notation $d_{\rm TV}$ to indicate the \textsl{total
variation distance} between laws of real valued random variables.

\begin{prop}[Lemma 5.3 in \cite{Chatterjee_ptrf}]\label{Chatterbox}
Assume that the random variable $W=f(Y)$ has zero mean and unit
variance, and denote by $Z \sim \mathscr{N}(0,1)$ a standard
Gaussian random variable. Then
$
d_{\rm TV}(W,Z)\leq 2 {\rm Var}(T(Y))^{1/2},
$
where the function $T(\cdot)$ is defined as
$$
T(y) = \int_0^1 \frac{1}{2\sqrt{t}}\sum_{i=1}^{n}E\big
[\frac{\partial f}{\partial y_i}(y)\frac{\partial f}{\partial
y_i}(\sqrt{t}\,y+\sqrt{1-t}\,\,Y)\big ]dt.
$$
\end{prop}
In what follows, we shall use Theorem \ref{theo:majDist} in order
to deduce a multidimensional generalization of Proposition
\ref{Chatterbox} (with the Wasserstein distance replacing total
variation).

\begin{prop}\label{MultiChatterbox}
Let $Y\sim \mathscr{N}_n(0, K)$, where $K=\{K(i,l):i,l=1,...,n\}$
is a $n\times n$ positive definite matrix. Consider absolutely
continuous functions $f_j : \real^n \rightarrow \real$,
$j=1,...,d$. Assume that each random variable $f_j(Y)$ has
zero mean, and also that each function $f_j$ and its partial
derivatives have subexponential growth at infinity. Denote by $Z
\sim \mathscr{N}_d(0,C)$ a Gaussian vector with values in
$\real^d$ and with positive definite covariance matrix $C=\{C(a,b)
: a,b=1,...,d\}$. Finally, write $W =(W_1,...,W_d) =
(f_1(Y),...,f_d(Y))$. Then,
$$
d_{\rm W}(W,Z)\leq \| C^{-1}\|_{op} \,\, \|
C\|_{op}^{1/2}\sqrt{\sum_{a,b=1}^d E[(C(a,b)-T_{ab}(Y))^2]}
$$
where the functions $T_{ab}(\cdot)$ are defined as
$$
T_{ab}(y) = \int_0^1 \frac{1}{2\sqrt{t}}\sum_{i,j=1}^{n}K(i,j)
E\big [\frac{\partial f_a}{\partial
y_i}(y)\frac{\partial f_b}{\partial
y_j}(\sqrt{t}\,y+\sqrt{1-t}\,\,Y)\big ]dt.
$$
\end{prop}
\begin{proof}
Without loss of generality, we can assume that
$Y=(Y_1,...,Y_n) = (X(h_1),...,X(h_n))$, where $X$ is an isonormal
Gaussian process over some Hilbert space $\HH$, and $\langle h_i ,
h_l\rangle_\HH = E(X(h_i)X(h_l))=K(i,l)$. According to Theorem
\ref{theo:majDist}, it is therefore sufficient to show that, for
every $a,b=1,...,d$,
$
T_{ab}(Y) = \langle DW_a , -DL^{-1}W_b \rangle_{\HH}.
$
To prove this last claim, introduce the two $\HH$-valued functions
$\Theta_a(y)$ and $\Theta_b(y)$, defined for $y\in\real^d$ as
follows:
$$
\Theta_a(y) =\sum_{i=1}^n \frac{\partial f_a}{\partial
y_i}(y)h_i
\quad\mbox{and}\quad
\Theta_b(y) = \int_0^1 \frac{1}{2\sqrt{t}}\sum_{j=1}^{n}\left\{
E\big [\frac{\partial f_b}{\partial
y_j}(\sqrt{t}\,y+\sqrt{1-t}\,\,Y)\big ]h_j \right\} dt.
$$
By using (\ref{baldassarre}), it is easily seen that $\Theta_a(Y)
= DW_a$. Moreover, by using e.g. formula (3.46) in
\cite{NourdinPeccati1}, one deduces that $\Theta_b(Y) =
-DL^{-1}W_b$. Since $T_{ab}(Y) =\langle \Theta_a(Y) , \Theta_b(Y)
\rangle_{\HH}$, the conclusion is immediately obtained.
\end{proof}

By specializing the previous statement to the case $n=d$ and
$f_j(y) = y_j$, $j=1,...,d$, one obtains the following simple
bound on the Wasserstein distance between Gaussian vectors of the
same dimension (the proof is straightforward and omitted).
\begin{corollary}
Let $Y\sim \mathscr{N}_d(0, K)$ and $Z\sim \mathscr{N}_d(0, C)$,
where $K$ and $C$ are two positive definite covariance matrices.
Then
$
d_{\rm W}(Y,Z)\leq  Q(C,K)\times \|C-K\|_{H.S.},
$
where
$$
Q(C,K) :=\min\{\| C^{-1}\|_{op} \,\, \| C\|_{op}^{1/2}, \|
K^{-1}\|_{op} \,\, \| K\|_{op}^{1/2}\}.
$$
\end{corollary}

\bigskip

\noindent \textbf{Acknowledgments.}
We thank an anonymous referee for insightful remarks.
Part of this paper has been
written while I. Nourdin and G. Peccati were visiting the
Institute of Mathematical Sciences of the National University of
Singapore, in the occasion of the ``Workshop on Stein's Method'',
from March 31 until April 4, 2008. These authors heartily thank A.
Barbour, L. Chen, K.-P. Choi and A. Xia for their kind hospitality
and support.


\begin{thebibliography}{99}

\bibitem{Barbour} A. D. Barbour (1990). Stein's method for diffusion approximations. \textit{Probab. Theory Related
Fields} \textbf{84}(3), 297--322.

\bibitem{BM} P. Breuer and P. Major (1983). \textrm{Central limit theorems
for nonlinear functionals of Gaussian fields.} \textit{J.
Multivariate Anal.} \textbf{13}(3), 425-441.

\bibitem{Chatterjee_ptrf} S. Chatterjee (2009). Fluctuation of eigenvalues and second order Poincar\'{e} inequalities.
\textit{Probab. Theory Related Fields} \textbf{143}, 1-40.

\bibitem{ChatterjeeMeckes} S. Chatterjee and E. Meckes (2008).
\textrm{Multivariate normal approximation using exchangeable
pairs.}
\textit{ALEA} {\bf 4}, 257-283.

\bibitem{Chen_Shao_sur} L. Chen and Q.-M. Shao (2005).
Stein's method for normal approximation. In: \textit{An
introduction to Stein's method}, 1-59. Lect. Notes Ser. Inst.
Math. Sci. Natl. Univ. Singap. \textbf{4}, Singapore Univ. Press,
Singapore, 2005.

\bibitem{DM} R.L. Dobrushin and P. Major (1979). \textrm {Non-central limit
theorems for nonlinear functionals of Gaussian fields.}
\textit{Z. Wahrsch. verw. Gebiete} \textbf{50}, 27-52.


\bibitem{Dudley book} R.M.\ Dudley (2003). \textit{Real Analysis and
Probability }(2$^{\text{nd}}$ Edition). Cambridge University
Press, Cambridge.

\bibitem{GS} L. Giraitis and D. Surgailis (1985). \textrm{CLT and other
limit theorems for functionals of Gaussian processes.}  \textit{Z.
Wahrsch. verw. Gebiete} \textbf{70}, 191-212.

\bibitem{Goetze} F. G\"{o}tze (1991). On the rate of convergence in the multivariate CLT.
\textit{Ann. Prob.} \textbf{19}(2), 724-739.

\bibitem{Hsu} E.P. Hsu (2005). Characterization of Brownian motion on manifolds
through integration by parts. In: \textit{Stein's method and
applications}, 195--208. Lect. Notes Ser. Inst. Math. Sci. Natl.
Univ. Singap. \textbf{5}. Singapore Univ. Press, Singapore.

\bibitem{NouPec07}
I. Nourdin and G. Peccati (2007a). Non-central convergence of
multiple integrals.
To appear in: \textit{Ann. Probab.}.

\bibitem{NourdinPeccati1} I. Nourdin and G. Peccati (2007b). \textrm{Stein's method on Wiener chaos.}
To appear in: \textit{Probab. Theory Related Fields}.

\bibitem{NourdinPeccati2} I. Nourdin and G. Peccati (2008).
\textrm{Stein's method and exact Berry-Ess\'een asymptotics for
functionals of Gaussian fields.} \textrm{Preprint. }

\bibitem{Nualart}
\rm D. Nualart (2006). {\it The Malliavin calculus and related
topics} (2$^{\text{nd}}$ Edition). \rm Springer Verlag, Berlin.

\bibitem{NualartOrtizLatorre} D. Nualart and S. Ortiz-Latorre (2008).
\textrm{Central limit theorem for multiple stochastic integrals
and Malliavin calculus}. \textit{Stochastic Process. Appl.} {\bf
118}(4), 614-628.

\bibitem{NP}
D. Nualart and G. Peccati (2005). Central limit theorems for
sequences of multiple stochastic integrals. \it Ann. Probab. {\bf
33} \rm (1), 177-193.

\bibitem{Pecc}
G. Peccati (2007). Gaussian approximations of multiple integrals.
\textit{Electron. Comm. Prob.} \textbf{12}, 350-364 (electronic).

\bibitem{PeccatiTudor} G. Peccati and C.A. Tudor (2005). Gaussian limits for vector-valued
multiple stochastic integrals. \textrm{ in}  \textit{S\'eminaire
de Probabilit\'es XXXVIII}, Lecture Notes in Math., Springer,
Berlin, 247--262.

\bibitem{Reinert}
\rm G. Reinert (2005). \rm  Three general approaches to Stein's
method. \rm in \textit{Lect. Notes Ser. Inst. Math. Sci. Natl.
Univ. Singap.}, \textbf{4}, 183--221. \rm Singapore Univ. Press,
Singapore.

\bibitem{ReinertRollin}
G. Reinert and A. R\"{o}llin (2007). Multivariate normal
approximation with Stein's method of exchangeable pairs under a
general linearity condition. Preprint available at
http://arxiv.org/abs/0711.1082

\bibitem{Stein1}
\rm Ch. Stein (1972). A bound for the error in the normal
approximation to the distribution of a sum of dependent random
variables. \rm In \textit{Proceedings of the Sixth Berkeley
Symposium on Mathematical Statistics and Probability (Univ.
California, Berkeley, Calif., 1970/1971), Vol. II: Probability
theory}, 583--602. \rm Univ. California Press, Berkeley, Calif..

\bibitem{Stein2}
\rm Ch. Stein (1986). {\it  Approximate computation of
expectations.} \rm Institute of Mathematical Statistics Lecture
Notes, Monograph Series, 7. \rm Institute of Mathematical
Statistics, Hayward, CA.


\end{thebibliography}
\end{document}